\documentclass[10pt]{article}
\usepackage{amsthm,amsfonts,amsmath,amscd,amssymb,epsfig}
\theoremstyle{plain}
\newtheorem{theorem}{Theorem}[section]
\newtheorem{thm}[theorem]{Theorem}

\newtheorem{lemma}[theorem]{Lemma}

\theoremstyle{remark}

\newtheorem{example}[theorem]{Example}

\theoremstyle{definition}

\newcommand{\id}{{{\mathchoice {\rm 1\mskip-4mu l} {\rm 1\mskip-4mu l}
{\rm 1\mskip-4.5mu l} {\rm 1\mskip-5mu l}}}}

\newcommand{\p}{\partial}

\newcommand{\into}{\hookrightarrow}

\newcommand{\Z}{{\mathbb{Z}}}
\newcommand{\R}{{\mathbb{R}}}
\newcommand{\C}{{\mathbb{C}}}

\newcommand{\tb}{{\rm tb}}
\title{Subcritical Stein manifolds are split}
\author{K.~Cieliebak}
\date{April 30, 2002}
\begin{document}
\maketitle
%\tableofcontents

%%%%%%%%%%%%%%%%%%%%%%%%%%%%%%%%%%%%%%%%%%%%%%%%%%%%%%%%%%%%%%%%%%%%%
\section{Introduction}
%%%%%%%%%%%%%%%%%%%%%%%%%%%%%%%%%%%%%%%%%%%%%%%%%%%%%%%%%%%%%%%%%%%%%

This note concerns the following little observation.

\begin{thm}\label{thm:def}
Every subcritical Stein manifold (domain) is deformation equivalent to
a split one.  
\end{thm}

Let me explain the words in the statement. 

{\em Stein. }A {\em Stein manifold} is a complex manifold $(W,J)$
without boundary
which admits an exhausting plurisubharmonic function
$\phi:W\to\R$. Here "exhausting" means proper and bounded from below,
and "plurisubharmonic" means that $-d(d\phi\circ J)$ is a positive
$(1,1)$ form. A sublevel set $\phi^{-1}(-\infty,c]$ in a Stein
manifold is called {\em Stein domain}. Note that Stein domains are
compact, whereas Stein manifolds are noncompact.  

{\em Subcritical. }It is well-known that the critical points of a
plurisubharmonic Morse function $\phi$
have Morse index $\leq n$, where $n$ is the complex dimension of $W$.
The Stein manifold (domain) is called {\em subcritical} if it admits
an exhausting plurisubharmonic Morse function 
$\phi$ all of whose critical points have index $<n$.

{\em Split. }A Stein manifold is called {\em split} if it is of the form
$(V\times\C,J_0\times i)$ for some Stein manifold $(V,J_0)$. An exhausting
plurisubharmonic function is given by $\phi=\phi_0+|z|^2$, where
$\phi_0$ is an exhausting plurisubharmonic function on $V$ and $z$ is the
complex coordinate on $\C$. A Stein domain is called split if it is a
sublevel set $\phi^{-1}(-\infty,c]$ in a split Stein manifold.
Note that split Stein manifolds (domains) are subcritical. 

{\em Deformation equivalence.} Two Stein structures $J_0,J_1$ on the
same smooth 
manifold $W$ are called {\em Stein homotopic} if there exists a
continuous family of Stein structures $J_t$ with exhausting
plurisubharmonic functions $\phi_t$ such that critical points of
$\phi_t$ do not travel to infinity during the homotopy. Two Stein
manifolds $(W,J)$ and 
$(W',J')$ are called {\em deformation equivalent} if there exists a
diffeomorphism $f:W\to W'$ such that $J$ and $f^*J'$ are Stein
homotopic on $W$. 
\medskip

A Stein manifold $(W,J,\phi)$ carries a canonical symplectic form
$-d(d\phi\circ J)$. Deformation equivalence implies
symplectomorphism~\cite{El2},
so this is the right notion of equivalence from the symplectic point
of view. Subcritical Stein manifolds have recently received some
interest because of their particularly simple symplectic properties
(\cite{BC},\cite{Mo},\cite{Vi},\cite{Ya}). In view of
Theorem~\ref{thm:def}, one can always assume 
them to be split. This will simplify the study of their symplectic
properties, e.g.~their symplectic field theory invariants~\cite{EGH}.  
\medskip

In the splitting $W=V\times\C$ the homeomorphism type of $V$ may not
be unique, as the following example shows.

\begin{example}\label{ex:nonunique}
Let $W=S^2\times\R^4$ with an almost complex structure $J$ of first
Chern class $c_1(J)=2k\in H^2(W;\Z)\cong\Z$. By the results
in~\cite{El2}, $W$ carries a unique Stein structure in the homotopy
class of $J$ for which the function 
$|x|^2$ on the $\R^4$-factor is plurisubharmonic. Then for every oriented
2-plane bundle $V\to S^2$ whose Euler class $e\in H^2(S^2;\Z)\cong\Z$ is
even and $e \leq -2-2k$, there exists a Stein structure $J_0$ on $V$
such that $V\times\C$ is deformation equivalent to $W$. 

To see this, consider first $V=T^*S^2$ with its natural Stein
structure $J_0$. It has Chern class $c_1(J_0)=0$ and Euler class (as a
bundle over $S^2$) $e=-2$. It is obtained from the 4-ball by attaching
a 2-handle along the Legendrian unknot with Thurston-Bennequin number
$\tb=-1$ and rotation number $r=0$ (see~\cite{GS}). Adding $l$
kinks to the knot, the Thurston-Bennequin number becomes $-1-l$, while
the rotation number $r$ can be made any of the integers
$-l,-l+2,\dots,l-2,l$. For $l$
even, the resulting Stein surface is the oriented 2-plane bundle over
$S^2$ with Euler class $e=\tb-1=-l-2$ and Chern class $c_1(J_0)=r$. So
for $e$ even and $e \leq -2k-2$ we can arrange $c_1(J_0)=2k$. By the
arguments below, $V\times\C$ is deformation equivalent to $W$. 

Note that the total spaces of 2-plane bundles over $S^2$ of different
Euler class are not 
homeomorphic. This can be seen, e.g., from the inverse limits of the
fundamental groups of the complements of compact subsets. So there are
infinitely many pairwise non-homeomorphic Stein surfaces $V$ such that
$V\times\C$ is Stein equivalent to $(W,J)$. 
\end{example}
\medskip

One could also ask whether a subcritical Stein manifold is
biholomorphic to a 
split one. Clearly a necessary condition is the existence of a
holomorphic embedding of $\C$ through every point. For example,
bounded open domains in $\C^n$ can never be biholomorphically
split. The following question grew out of a discussion with R.~Hind. 

{\em Is every subcritical Stein manifold which has a holomorphic
embedding of $\C$ through every point biholomorphic to a split one?}

{\bf Acknowledgement. }The question whether every subcritical Stein
manifold is split was brought up by K.~Mohnke. This note was written
during my stay at the Institute for Advanced Study, which I thank for
its hospitality.

%%%%%%%%%%%%%%%%%%%%%%%%%%%%%%%%%%%%%%%%%%%%%%%%%%%%%%%%%%%%%%%%%%%%%
\section{Proof}
%%%%%%%%%%%%%%%%%%%%%%%%%%%%%%%%%%%%%%%%%%%%%%%%%%%%%%%%%%%%%%%%%%%%%

The proof is based on Eliashberg's theory of Stein manifolds
(\cite{El1},\cite{El2}). Recall first some notation from~\cite{El1}. Let
$(W,J)$ be a Stein domain of complex dimension $n$. A {\em
handle attaching triple (HAT)} $(f,\beta,\gamma)$ is the data for
attaching a handle in the catagory of almost complex manifolds. Here
$f:S^{k-1}\into\p W$ is an embedding. It induces 
injective bundle homomorphisms $df:TS^{k-1}\to f^*T(\p W)$ and 
$Df:TS^{k-1}\oplus\underline{\R}\to f^*TW$ by sending the generator of
the trivial bundle $\underline{\R}$ to an inward pointing vector field
transverse to $\p W$. $\beta$ is a normal framing for
$f$ in $\p W$, i.e.~a bundle isomorphism
$\beta:S^{k-1}\times\R^{2n-k}\to\nu_f$ to the normal bundle $\nu_f$ to
$f$ in $\p W$. $\gamma:S^{k-1}\times\C^n\to f^*TW$ is an isomorphism
of complex bundles which is homotopic to $Df\oplus\beta$ as an
isomorphism of real bundles. Here
$TS^{k-1}\oplus\underline{\R}\oplus\R^{2n-k}$ is identified with
$S^{k-1}\times\C^n$ 
by viewing $TS^{k-1}\oplus\underline{\R}$ as the tangent bundle to the
ball $B^k\times 0\times 0\subset \R^k\times i\R^k\times\C^{n-k}$.  

An {\em isotopy} of HATs is an isotopy of embeddings $f_t$ covered by
homotopies $\beta_t,\gamma_t$. Attaching a handle with isotopic HATs
yields diffemorphic smooth manifolds with homotopic almost complex
structures. 

Note that the homotopy group $\pi_{k-1}SO(2n-k)$ acts
transitively on the framings $\beta$ (considered up to homotopy) by
composition. For a HAT $(f,\beta,\gamma)$ and $g$ in the kernel of the
map $\pi_{k-1}SO(2n-k) \to \pi_{k-1}SO(2n)$, $(f,\beta\cdot
g,\gamma)$ is again a HAT (which leads to a different smooth manifold
when attaching the handle). 

The maximal complex subspaces of $\p W$ define a contact structure
$\xi$ on $\p W$. An embedding into $\p W$ is called {\em isotropic} if
it is tangent to $\xi$. An isotropic embedding of the maximal possible
dimension, $n-1$, is called {\em Legendrian}. A HAT $(f,\beta,\gamma)$
is called {\em special} if $f$ is an isotropic embedding, 
$\beta=JDf\oplus\theta$ for an 
injective complex bundle homomorphism $\theta:S^{k-1}\times\C^{n-k}\to
f^*TW$, and $\gamma=Df\oplus JDf\oplus\theta$. It is called {\em stably
special} if there exists a $g\in \ker[\pi_{k-1}SO(2n-k) \to
\pi_{k-1}SO(2n)]$ such that $(f,\beta\cdot
g,\gamma)$ is special. The main inductional lemma in~\cite{El1} states
that if $(W,J)$ is a Stein domain and $(f,\beta,\gamma)$ a special
HAT, then a handle can be attached in such a way that the Stein
structure extends over the handle.   

Now let $(V,J_0)$ be a Stein domain of complex dimension $n-1$ with
plurisubharmonic Morse function $\phi_0$, $\phi_0|_{\p V}\equiv c$. Let
$W\subset V\times\C$ be the Stein domain
$(\phi_0+|z|^2)^{-1}(-\infty,c]$, where $z$ is the coordinate on
$\C$. $W$ is equipped with the complex structure $J=J_0\times i$ and
the plurisubharmonic function $\phi=\phi_0+|z|^2$. 
Note that $\p W$ has a natural open book structure with trivial
monodromy,
$$\p W\cong \p V\times B^2\cup V\times S^1.$$ 
Let $(f,\beta,\gamma)$ be a HAT for $W$ of index $k<n$. 

A HAT $(f_0,\beta_0,\gamma_0)$ for $V$ naturally induces a HAT 
$(\hat f_0,\hat \beta_0,\hat \gamma_0)$ for $W$ with $\hat
f_0=f_0\times 0:S^{k-1}\into\p V\times 0$,
$\hat\beta_0=\beta_0\times\id_\C$ and
$\hat\gamma_0=\gamma_0\times\id_\C$.

\begin{lemma}\label{lem:HAT1}
There exists an embedding $f_0:S^{k-1}\into\p V$ such that 
$\hat f_0$ is isotopic (through embeddings into $\p W$) to $f$. 
\end{lemma}

\begin{proof}
Let $\Delta\subset V$ be the skeleton, i.e.~the union of all
descending manifolds of critical points of $\phi_0$ (with respect to
some Riemannian metric). Note that the
negative gradient flow of $\phi_0$ retracts $V$ onto $\Delta$. Since
$\Delta\times S^1$ has codimension at least $n-1$ and $S^{k-1}$ has
dimension $<n-1$, a small perturbation of $f:S^{k-1}\into\p W$ makes
it avoid 
$\Delta\times S^1$. Then we use the gradient flow of $\phi_0$ to isotop
$f$ to an embedding into $\p V\times B^2$. Since $2\dim
S^{k-1}\leq\dim\p V-1$, the latter embedding is isotopic to an
embedding into $\p V\times 0$.
\end{proof}

\begin{lemma}\label{lem:HAT2}
There exist a HAT $(f_0,\beta_0,\gamma_0)$ on $V$
such that the HAT $(\hat f_0,\hat \beta_0,\hat \gamma_0)$ is
isotopic to $(f,\beta,\gamma)$.
\end{lemma}

\begin{proof}
The previous lemma shows that after an isotopy of HATs, we may assume
that $f=\hat f_0$. Complete $f_0$ to any HAT
$(f_0,\bar\beta_0,\bar\gamma_0)$ on $V$. This allows us to identify
homotopy classes of framings $\beta_0$ for $f_0$ with
$\pi_{k-1}SO(2n-k-2)$, and similarly for
$\gamma_0,\beta,\gamma$. Now consider the following commutative
diagram: 
\begin{equation*}
\begin{CD}
        \beta_0\in\pi_{k-1}SO(2n-k-2) @>\sigma>>
        \pi_{k-1}SO(2n-k)\ni\beta \\ 
        @VVV @VV{\cong}V \\
        \pi_{k-1}SO(2n-2) @>{\cong}>> \pi_{k-1}SO(2n) \\ 
        @AAA @AAA \\
        \gamma_0\in\pi_{k-1}U(n-1) @>{\cong}>> \pi_{k-1}U(n)\ni\gamma.  
\end{CD}
\end{equation*}
Here the map $\sigma$ is surjective by Bott periodicity
(see~\cite{Ko}, p.~230): $\pi_iO(l-1)\to\pi_iO(l)$ is an isomorphism
for $i<l-2$ and surjective for $i\leq l-2$. With 
$i=k-1$ and $l-1=2n-k-2$, the condition for surjectivity becomes
$k-1\leq 2n-k-3$, or $k\leq n-1$, which is fulfilled by hypothesis.
The same argument yields the isomorphisms in the middle row and at the
vertical arrow. The isomorphism in
the bottom row follows similarly (see~\cite{Mi}):
$\pi_iU(n-1)\to\pi_iU(n)$ is 
surjective for $i\leq 2n-2$ and an isomorphism for $i<2n-2$. 

We see that given $\beta,\gamma$ with the same image in
$\pi_{k-1}SO(2n)$ we find preimages $\beta_0,\gamma_0$ under the
vertical maps. Since $\beta_0,\gamma_0$ have the same image in
$\pi_{k-1}SO(2n-2)$, $(f_0,\beta_0,\gamma_0)$ is a HAT with the
desired properties. 

Note that here we make a choice if $\sigma$ is not bijective: We may
change $\beta_0$ by any element in the kernel of $\sigma$ and still
get a HAT with the desired properties. This freedom will be important
later on. 
\end{proof}

\begin{lemma}\label{lem:HAT3}
There exist a special HAT $(f_0,\beta_0,\gamma_0)$ on $V$
such that the HAT $(\hat f_0,\hat \beta_0,\hat \gamma_0)$ is
isotopic to $(f,\beta,\gamma)$.
\end{lemma}

\begin{proof}
Let $(f_0,\beta_0,\gamma_0)$ be a HAT as provided by
Lemma~\ref{lem:HAT2}. It is shown in~\cite{El1} that
$(f_0,\beta_0,\gamma_0)$ is isotopic (on $V$) to a stably special HAT
$(f_1,\beta_1,\gamma_1)$. Moreover, for $n>3$ or $k<n-1$, the stably
special HAT $(f_1,\beta_1,\gamma_1)$ is isotopic to a special one. The
same holds for $n=2$ by elementary reasons. 

It remains to treat the case $n=3$ and $k=2$. Then the diagram in the
proof of Lemma~\ref{lem:HAT2} becomes
\begin{equation*}
\begin{CD}
        \beta_1\in\pi_1SO(2) @>\sigma>>
        \pi_1SO(4)\ni\beta \\ 
        @VVV @VV{\cong}V \\
        \pi_1SO(4) @>{\cong}>> \pi_1SO(6) \\ 
        @AAA @AAA \\
        \gamma_1\in\pi_1U(2) @>{\cong}>> \pi_1U(3)\ni\gamma.  
\end{CD}
\end{equation*}
The embedding $f_1$ is a Legendrian knot in the contact 3-manifold $\p
V$. It determines a normal framing $JDf$ in $\p V$ which corresponds
to an element in $\pi_1SO(2)\cong\Z$. The HAT
is special iff $\beta_1=JDf\in\pi_1SO(2)\cong\Z$. Adding ``kinks'' to $f_1$
allows us to {\em decrease} $JDf$ by an arbitrary integer, 
but by Bennequin's inequality~\cite{Be}, we cannot increase
it. However, recall 
that in the proof of Lemma~\ref{lem:HAT2} we had the freedom of
changing $\beta_0$ (and hence $\beta_1$) by an arbitrary element in
$\ker[\sigma:\pi_1SO(2)\to\pi_1SO(4)]\cong 2\Z$. So we can make $\beta_1$
smaller than $JDf$ by subtracting an even integer and then decrease
$JDf$ by adding kinks to $f_1$ until $JDf=\beta_1$. The result is a
special HAT as desired.  

Note that the modification of $\beta_1$ changes the diffeomorphism
type of the Stein surface we get after attaching the handle. Since we
have the freedom of making $\beta_1$ more negative and adding more
kinks to $f_1$, the diffeomorphism type is not uniquely determined
(see Example~\ref{ex:nonunique}).  
\end{proof}

\begin{proof}[Proof of Theorem~\ref{thm:def}]
The main extension lemma in~\cite{El1} states that given a special
HAT $(f_0,\beta_0,\gamma_0)$, a handle can be attached in such a way
that the Stein structure extends over the handle. So using
Lemma~\ref{lem:HAT3} and induction, we obtain a Stein manifold
$(V,J_0)$ and a diffeomorphism $F:V\times\C\to W$ such that $F^*J$ is
homotopic (through almost complex structures) to $J_0\times
i$. Moreover, by construction, there are exhausting plurisubharmonic
functions $\phi$ on $W$ and $\phi_0$ on $V$ such that $F^*\phi$ and
$\phi_0+|z|^2$ have the same critical points. Now another theorem of
Eliashberg~\cite{El2} implies that $J_0\times i$ and $f^*J$ are Stein
homotopic. 

The same arguments work for Stein domains. This proves
Theorem~\ref{thm:def}. 
\end{proof}

%%%%%%%%%%%%%%%%%%%%%%%%%%%%%%%%%%%%%%%%%%%%%%%%%%%%%%%%%%%%%%%%%%%%%%%
%%%%%%%%%%%%%%%%%%%%%%%%%%%%% References %%%%%%%%%%%%%%%%%%%%%%%%%%%%%%
%%%%%%%%%%%%%%%%%%%%%%%%%%%%%%%%%%%%%%%%%%%%%%%%%%%%%%%%%%%%%%%%%%%%%%%

\end{document}